%% file: nuf.tex
\input btxmac.tex

\input nuf.frm

\line{\ \hfill \ }
\bigskip

\centerline{\titlefont Non-unique factorization of polynomials}
\smallskip
\centerline{\titlefont over residue class rings of the integers}
\bigskip

\centerline{\authorfont Christopher Frei\ \ and\ \ Sophie Frisch}
\bigskip\bigskip

\input nuf.bdy
\bigskip

\subhead{References}%
\bibliography{nuf}
\bibliographystyle{siamese}

\bigskip\bigskip
\goodbreak
{
\obeylines\frenchspacing
Institut f\"ur Mathematik A
Technische Universit\"at Graz
Steyrergasse 30
A-8010 Graz, Austria
\smallskip
Phone +43 316 873 7133
Fax +43 316 873 107133
\bigskip
\tt
\line{\vbox{\hbox{Ch.~F.}\hbox{\tt frei@math.tugraz.at}}\hfill
\vbox{\hbox{S.~F.}\hbox{\tt frisch@blah.math.tu-graz.ac.at}} }
}

\bye

%% file: nuf.bdy
{\csc Abstract.}
We investigate non-unique factorization of polynomials in
${\bb Z}_{p^n}[x]$ into irreducibles. As a Noetherian ring whose
zero-divisors are contained in the Jacobson radical, ${\bb Z}_{p^n}[x]$
is atomic. We reduce the question of factoring arbitrary non-zero
polynomials into irreducibles to the problem of factoring monic polynomials
into monic irreducibles.
The multiplicative monoid of monic polynomials of ${\bb Z}_{p^n}[x]$ is a
direct sum of monoids corresponding to irreducible polynomials in
${\bb Z}_{p}[x]$, and we show that each of these monoids has infinite
elasticity.
Moreover, for every $m\in {\bb N}$, there exists in each of these monoids
a product of $2$ irreducibles that can also be represented as a product of
$m$ irreducibles. 


\medskip
\footnote{}{Author Posting. (c) Taylor \& Francis Group, LLC, 2011.}
\footnote{}{This is the author's version of the work. It is posted here by permission of Taylor \& Francis Group, LLC,
for personal use, not for redistribution.}
\footnote{}{The definitive version was published in Communications in Algebra, Volume 39 Issue 4, April 2011.
doi:10.1080/00927872.2010.549158}
\subhead{1. Introduction}

Polynomials with coefficients in $\Zpn=\intz/p^n\intz$, $p$ prime,
$n\ge 2$, have been investigated quite a lot with respect to the
polynomial functions they define \cite{JPSZ06PFFCR,Zh04PFPP,Fri01PFFCR}.
Comparatively little is known about the multiplicative monoid of
non-zerodivisors of $\Zpnx$.
McDonald \cite{McD74FRI} and Wan \cite{Wan03FFGR} list a few properties
concerning uniqueness of certain kinds of factorizations of polynomials in
$\Zpnx$. At the same time, factorization into irreducibles is manifestly
non-unique. In this paper, we investigate non-uniqueness of factorization
of polynomials in $\Zpnx$ into irreducibles, in the spirit of \cite{HKG06NUF}.
Throughout this paper, $p$ is a fixed prime and $n\ge 2$.

Here is an example of the phenomena we will study. Consider the
identity in $\Zpnx$:
$$(x^m+p^{n-1})(x^m+p) = x^m (x^m + p^{n-1} + p).$$
Assuming for the moment that every element of $\Zpnx$ has a factorization
into irreducible elements (cf.~\citeZpnxatomic), we note that
$(x^m+p^{n-1})$ is a product of at most $n-1$ irreducibles, and that
$(x^m+p)$ is irreducible.
This is so because both polynomials represent a power of $x$ in $\Zpx$. 
By unique factorization in $\Zpx$, each of their factors in $\Zpnx$
(apart from units) also represents a power of $x$ in $\Zpx$. 
(We have used the fact that a polynomial in $\Zpnx$ is a unit if and
only if it maps to a unit in $\Zpx$ under canonical projection.) 
Since every non-unit factor in $\Zpnx$ of $(x^m+p^{n-1})$ or 
$(x^m+p)$ represents a power of $x$ in $\Zpx$, the constant coefficient
of every such factor is divisible by $p$.
On the other hand,
the constant coefficient of $(x^m+p^{n-1})$ is divisible by no higher
power of $p$ than $p^{n-1}$, so it cannot be a product of more than
$n-1$ non-unit polynomials, and $(x^m+p)$ cannot be a product of more
than one non-unit in $\Zpnx$.
We have seen that for arbitrary $m\in\natn$, there exists in $\Zpnx$
a product of at most $n$ irreducibles that is also representable
as a product of more than $m$ irreducibles.

We now establish some notation:

\rmprofess{\vDef} 
We denote $p$-adic valuation on $\ratq$ by $v$ and extend 
$v$ to $\ratq(x)$ by defining $v(f)=\min_k v(a_k)$ for
$f=\sum_k a_kx^k\in\intz[x]$, and $v(f/g)=v(f)-v(g)$.
This gives a surjective mapping 
$v\colon \intz[x]\rightarrow \natn_0\cup\{\infty\}$.

We denote by $(N_n,+,\le)$ the ordered monoid with elements
$0$, $1,\ldots$, ${n-1}$, $\infty$ resulting from factoring 
$(\natn_0\cup\{\infty\},+,\le)$ by the congruence relation
that identifies all values greater or equal $n$, including $\infty$.
By abuse of notation, we still use $v$ for the surjective
mapping $v\colon \Zpnx\rightarrow N_n$ obtained by factoring
$p$-adic valuation $v\colon \intz[x]\rightarrow \natn_0\cup\{\infty\}$
by the above congruence relation.
\endrmprofess

The mapping $v\colon \Zpnx\rightarrow N_n$ behaves like a valuation,
except that $(N_n,+,\le)$ is not a group (nor can it be extended to
a group, as it is not cancellative). 

\profess{\vproperties}
$v\colon \Zpnx\rightarrow N_n$ satisfies
{\itemenv
\item{\rm (i)} $v(f)=\infty \iff f=0$
\item{\rm (ii)} $v(f+g)\ge \min(v(f), v(g))$
\item{\rm (iii)} $v(f\cdot g)=v(f)+ v(g)$
\par
}
\endprofess

We will write $p$ for a prime in $\intz$ as well as for its residue
class in $\Zpn$ or $\Zpnx$. 
This should not cause any confusion; we are always going to specify
which ring we are talking about. In the remainder of this section,
we review a few well-known facts about $\Zpnx$.

\profess{\reviewZpnx}
{\itemenv
For $f\in\Zpnx$, the following are equivalent:
\itemitem{\rm(a)} $v(f)>0$, i.e., all coefficients of $f$ are
divisible \(in $\Zpn$\) by $p$.
\itemitem{\rm(b)} $f$ is nilpotent.
\itemitem{\rm(c)} $f$ is a zero-divisor.
\par
}
\endprofess

\proof
This follows from the properties of $v \colon \Zpnx\rightarrow N_n$
mentioned in \vpropertiesNum.
\endproof

\rmprofess{\NilJZDef} Let $R$ be a commutative ring.
{\itemenv
\item{(i)} $\NilR$ denotes the nilradical of $R$, i.e., the intersection
of all prime ideals of $R$.
\item{(ii)} $\J(R)$ denotes the Jacobson radical of $R$, i.e., the
intersection of all maximal ideals of $R$.
\item{(iii)} $\Z(R)$ denotes the set of zero-divisors of $R$.
\item{(iv)} $\U(R)$ denotes the group of units of $R$.
\par
}
\endrmprofess

As we all know, the intersection of all prime ideals of a 
commutative ring coincides with the set of nilpotent elements:
$$\NilR\kern6pt = \kern-6pt \bigcap_{P \hbox{\rm\ prime}} \kern-8pt P 
\kern6pt =
\kern6pt \{r\in R\mid \exists n\in\natn\;\; r^n=0\}.$$

\profess{\zerodivnilJ}
$$\Nil(\Zpnx) = \J(\Zpnx) = (p) = \Z(\Zpnx)$$
\endprofess

\proof
Since the maximal ideals of $\Zpnx$ are precisely the ideals $(p,f)$ with
$f$ representing an irreducible polynomial in $\Zpx$, we see $\J(\Zpnx)=(p)$.
We already know $(p) = \Nil(\Zpnx) =\Z(\Zpnx)$ from \citereviewZpnx.
\endproof

\profess{\units}
Let $R$ be a commutative ring. The units of $R[x]$ are precisely
the polynomials $a_0+a_1x+\ldots+a_nx^n$ with $a_0$ a unit of $R$ and
$a_k$ nilpotent for $k>0$. 
\endprofess

\proof
If $f$ is a unit then its image under projection to $(R/P)[x]$ is
a unit for every prime ideal $P$ of $R$. Therefore, $a_0$ is not
in any $P$ and hence a unit; and for $k>0$, $a_k$ is in every
$P$ and therefore nilpotent. Conversely, if $f=a_0+g$ with $a_0$
a unit of $R$ and all coefficients of $g$ in the intersection of
all prime ideals of $R$, then $g$ is in every prime ideal of
$R[x]$ and hence $f=a_0+g$ is in no prime ideal of $R[x]$ and
therefore a unit of $R[x]$.
\endproof

\profess{\unitscor}
Let $p$ be a prime and $n\ge 1$.  Then
$\U(\Zpnx)$ consists of precisely the polynomials
$f=a_0+a_1x+\ldots+a_nx^n$ such that \(in $\Zpn$\)
$p\dividesnot a_0$ and $p\divides a_k$ for $k>0$.
It follows that a polynomial in $\intz[x]$ represents a
unit in $\Zpnx$ for some $n\ge 1$ if and only if it
represents a unit in $\Zpnx$ for all $n$.
\endprofess

\subhead{2. Different concepts of irreducibility in rings with zero-divisors}

Various facts about irreducible and prime elements that are clear
in integral domains tend to become messy and obfuscated as soon as
zero-divisors are involved. We will see, however, that everything 
works just as smoothly as in integral domains if all zero-divisors
are contained in the Jacobson radical. For an in-depth investigation
of factorization in rings with zero-divisors see 
\cite{AndVL96FZD,AndVL99FZDii,AgAnVL01FZDiii}. 
The very natural concepts defined in this section appear in the
literature (\cite{AndVL96FZD,AndVL99FZDii,AgAnVL01FZDiii} and their
references) under varying names. 
Our choice of names has been motivated by a desire to define
{\em associated} and {\em irreducible} for commutative rings just like 
they are usually defined for integral domains and to reflect the
existence of a group of ``stronger'' concepts (associated, irreducible,
atomic), and a parallel group of ``weaker'' concepts, with analogous
relationships between the concepts inside each group.
$R$ is always a commutative ring with identity.

\rmprofess{\IrredPrimeDef} 
If $c\in R$ is a non-zero non-unit then we say that $c$ is
{\em weakly irreducible} if for all $a,b\in R$
$$ c=ab\Longrightarrow c\divides a \hbox{\rm\ or\ } c\divides b$$
and that $c$ is {\em irreducible} if for all $a,b\in R$
$$ c=ab\Longrightarrow \hbox{\rm $b$ is a unit or $a$ is a unit}.$$
A non-zero non-unit $c$ is called {\em prime} if for all $a,b\in R$
$$ c\divides ab\Longrightarrow c\divides a \hbox{\rm\ or\ } c\divides b$$
\endrmprofess

\rmprofess{\AssocDef}
Let $a,b\in R$.
We call $a$ and $b$ {\em weakly associated} if $a\divides b$ and
$b\divides a$ (or equivalently, if $(a)=(b)$). 
We call $a$ and $b$ {\em associated} if there exists a unit 
$u\in R$ such that $a=bu$.
\endrmprofess

It is clear that irreducible implies weakly irreducible, that associated 
implies weakly associated, and that both converses hold in integral
domains. Also, prime implies irreducible in integral domains.
We now consider a class of commutative rings properly containing the
class of integral domains in which the implications mentioned above
still hold.

\rmprofess{\HarmlessDef} 
Let $R$ be a commutative ring. We say that $R$ is a ring with
{\em harmless zero-divisors}, if 
$\Z(R)\subseteq 1-\U(R)=\{1-u\mid \hbox{\rm $u$ a unit of $R$}\}$.
(These are the
rings called {\it pr\'esimplifiable} by Bouvier (cf.~\cite{AndVL96FZD}).)
\endrmprofess

\profess{\weakstrongequiv} {\rm(\cite{AndVL96FZD})}.
Let $R$ be a ring with harmless zero-divisors and $a,b,c,u,v\in R$.
Then
{\itemenv
\item{\rm (1)} if $a\ne 0$, $a=bu$ and $b=av$ then $u,v$ are units;
\item{\rm (2)} $a,b$ are weakly associated if and only if they are associated;
\item{\rm (3)} $c$ is weakly irreducible if and only if $c$ is irreducible;
\item{\rm (4)} if $c$ is prime then $c$ is irreducible.
\par
}
\endprofess

\proof
(1) $a=avu$; therefore $a(1-vu)=0$, which makes
$(1-vu)$ a zero-divisor. Now $1-vu=1-w$ for some unit $w$; hence
$vu=w$ and $u,v$ are units.
(2), (3) and (4) follow from (1).
\endproof

\rmprofess{\JcontainsZharmless}
If a commutative ring $R$ satisfies $\Z(R)\subseteq \J(R)$ then the 
statements of \citeweakstrongequiv\ hold for $R$. In particular,
they hold for $\Zpnx$.
\endrmprofess

\proof
It is easy to see that $\J(R)\subseteq 1-\U(R)$ holds in every
commutative ring. Therefore every ring satisfying $\Z(R)\subseteq \J(R)$
is a ring with harmless zero-divisors.
In particular, $\Zpnx$ is a ring with harmless zero-divisors, since
$\Z(\Zpnx) = \J(\Zpnx)$  by \citezerodivnilJ.
\endproof

\rmprofess{\AtomicDef}
Let $R$ be a commutative ring. $R$ is called {\em atomic} if every
non-zero non-unit is a product of irreducible elements. We call $R$
{\em weakly atomic} if every non-zero non-unit is a product of
weakly irreducible elements.
\endrmprofess

If one examines the usual proof that an integral domain with
ACC for principal ideals is atomic, one finds that it shows
that any commutative ring with ACC for principal ideals is weakly
atomic. Namely, suppose there exists an element $s$ in the set $S$
of those non-zero non-units of $R$ that do not factor into weakly
irreducible elements. 
Then $s$ is not weakly irreducible, so there exist $a,b$ with $s=ab$
and $s\dividesnot a$, $s\dividesnot b$.
Such $a,b$ must be non-zero non-units, and at least one of them is in $S$;
say $a\in S$. As $a\divides s$ and $s\dividesnot a$, we have
$(s)\subset (a)$. By iteration we get an infinite ascending chain
of principal ideals. We have shown:

\profess{\ACCweaklyatomic} {\rm(\cite{AndVL96FZD} Thm.~3.2)}.
If $R$ is a commutative ring satisfying ACC for principal ideals
then $R$ is weakly atomic.
\endprofess

Combining Lemmas \weakstrongequivNum\ and \ACCweaklyatomicNum, we obtain

\profess{\ACCatomic} {\rm (\cite{AndVL96FZD})}.
Every commutative ring with harmless zero-divisors satisfying ACC for
principal ideals is atomic; in particular, every Noetherian ring with
harmless zero-divisors is atomic.
\endprofess

In particular, \citeACCatomic\ applies to $\Zpnx$:

\profess{\Zpnxatomic}
$\Zpnx$ is atomic.
\endprofess

\subhead{3. Uniqueness of factorization}

In the previous section we have seen that every non-zero non-unit of
$\Zpnx$ is a product of irreducibles. Step by step we now reduce the task
of factoring arbitrary elements of $\Zpnx$: first to factoring 
non-zerodivisors, then to factoring monic polynomials and,
finally, to factoring monic primary polynomials. This is a matter
of applying well-known facts of the ``uniqueness'' kind from \cite{McD74FRI}
or \cite{Wan03FFGR}. It is only in the next section, concerned with
factoring monic primary polynomials, that non-uniqueness of
factorization comes into play.

\subsubhead{Arbitrary polynomials to non-zerodivisors.}


\profess{\onlyprimep}
Let $p$ be a prime and $n\ge 2$. 
For $f\in\Zpnx$, the following are equivalent:
{\itemenv
\item{\rm (1)} $f=pu$ for some unit $u$ of $\Zpnx$.
\item{\rm (2)} $f$ is prime.
\item{\rm (3)} $f$ is irreducible and a zero-divisor.
\par
}
\endprofess

\proof
$(1\Rightarrow 2)$.
Let $v\colon \Zpnx\rightarrow N_n$ be $p$-adic valuation as
defined in \vDefNum. Since $v(p)=1$, additivity of $v$
(\vpropertiesNum) implies that $p$ is prime.
Therefore every element associated to $p$ is prime.

$(2\Rightarrow 3)$.
By \citeJcontainsZharmless, prime elements of $\Zpnx$ are irreducible.
If $c$ is prime, then $(c)\supseteq {\rm Nil}(\Zpnx)=(p)$; so $c\divides p$
and therefore $c$ and $p$ are associated.
As $p$ is a zero-divisor, so is $c$.

$(3\Rightarrow 1)$.
If $c$ is a zero-divisor, then $c\in\Z(\Zpnx)=(p)$, meaning
$c=pv$ for some $v$. If, furthermore, $c$ is irreducible, then
$v$ must be a unit.
\endproof

\rmprofess{\cancelp}
Let $f\in \Zpnx$, $f\ne 0$. 
{\itemenv
\item{\rm (i)}
There exists a non-zerodivisor $g\in\Zpnx$, and $0\le k<n$, such that
$f=p^kg$. Furthermore, $k$ is uniquely determined by $k=v(f)$,
and $g$ is unique modulo $(p^{n-k})$. 
\item{\rm (ii)}
In every factorization of $f$ into irreducibles,
exactly $v(f)$ of the irreducible factors are associates of $p$.
\par
}
\endrmprofess

\proof
(i) follows from \citereviewZpnx\ and the definition of $p$-adic valuation.
(ii) follows from (i) and the fact that $p$ is prime in $\Zpnx$ 
(by \citeonlyprimep).
\endproof

\subsubhead{Nonzerodivisors to monic polynomials.}

\rmprofess{\monicrepresentative} {\rm(\cite{McD74FRI}, Thm.~XIII.6)}.
Every non-zerodivisor $f\in\Zpnx$ is uniquely representable as
$f=ug$, with $u$ a unit of $\Zpnx$ and $g$ monic. The degree of
$g$ is $\deg \bar f$, where $\bar f$ is the image of $f$ in
$\Zpx$ under canonical projection.
\endrmprofess

Actually, McDonald \cite{McD74FRI} only shows the existence of
$g$ and $u$ in \monicrepresentativeNum, but uniqueness is easy:
suppose $ug=vh$ with $u,v$ units, $g,h$ monic. Then $v^{-1}ug=h$.
As $g,h$ are monic, so is $v^{-1}u$. By \unitscorNum, the only
monic unit of $\Zpnx$ is $1$, however, so $u=v$ and $g=h$.

Applying \monicrepresentativeNum\ to the irreducible factors of a 
non-zerodivisor of $\Zpnx$, we obtain:

\rmprofess{\factoringregularf}
Let $f\in\Zpnx$, not a zero-divisor, and $u$ and $g$ the unique unit
and monic polynomial, respectively, in $\Zpnx$ with $f=ug$. 
For every factorization $f=c_1\cdot\ldots\cdot c_k$ into irreducibles,
there exist uniquely determined monic irreducible $d_1,\ldots, d_k\in\Zpnx$
and units $v_1,\ldots,v_k$ in $\Zpnx$ with $c_i=v_id_i$,
$u=v_1\cdot\ldots\cdot v_k$ and $g=d_1\cdot\ldots\cdot d_k$.
\endrmprofess

By \citefactoringregularf, we have reduced the 
question of factoring non-zerodivisors of $\Zpnx$ into irreducibles to
the question of factoring monic polynomials into monic irreducibles.
As there are only finitely many monic polynomials of each degree
in $\Zpnx$, we also see that every non-zerodivisor of $\Zpnx$ has
(up to associates) only finitely many factorizations into irreducibles.

\subsubhead{Monic polynomials to primary monic polynomials.}

\rmprofess{\primary}
Let $f\in \Zpnx$, not a zero-divisor. Then $(f)$ is a primary ideal
of $\Zpnx$ if and only if the image of $f$ under the canonical projection
$\pi\colon\Zpnx\rightarrow\Zpx$ is associated to a power of an
irreducible polynomial.
\endrmprofess

\proof
In the principal ideal domain $\Zpx$, the non-trivial primary ideals
are precisely the principal ideals generated by powers of irreducible
elements. The canonical projection $\pi$ induces a bijective 
correspondence between primary ideals of $\Zpx$ and primary ideals
of $\Zpnx$ containing $(p)$. 

Now an ideal of $\Zpnx$ that does not consist only of zero-divisors
is primary if and only if its radical is a maximal ideal (since the
only non-maximal prime ideal is $(p)=\Z(\Zpnx)$).
Let $f\in \Zpnx$. Then, since every prime ideal of $\Zpnx$ contains
$(p)$, the radical of $(f)$ is equal to the radical of
$(f)+(p)=\pi^{-1}(\pi((f)))$.

For a non-zerodivisor $f$, therefore, $(f)$ is primary if and only if
$(f)+(p)$ is primary, which is equivalent to $\pi(f)$ being a primary
element of $\Zpx$.
\endproof

\rmprofess{\PrimaryDef}
We call a non-zerodivisor of $\Zpnx$ primary if its image under
projection to $\Zpx$ is associated to a power of an irreducible
polynomial.
\endrmprofess

The simplest form of Hensel's Lemma shows that every monic
$f\in \Zpnx$ is a product of primary polynomials.
Furthermore, the monic primary factors of a monic polynomial in $\Zpnx$
are uniquely determined.

\rmprofess{\uniqueprimary} {\rm(\cite{Wan03FFGR}, Thm.~13.8)}. 
Let $f\in \Zpnx$ monic, then there exist monic polynomials
$f_1,\ldots, f_r \in \Zpnx$ such that $f=f_1\cdot\ldots\cdot f_r$
and the residue class of $f_i$ in $\Zpx$ is a power of
a monic irreducible polynomial $g_i\in \Zpx$, with
$g_1,\ldots, g_r$ distinct.
The polynomials $f_1,\ldots, f_r \in \Zpnx$ are \(up to ordering\)
uniquely determined.
\endrmprofess

\profess{\uniquefactorization}
Every non-zero polynomial $f\in \Zpnx$ is representable as
$$f=p^kuf_1\cdot\ldots\cdot f_r$$
with $0\le k<n$, $u$ a unit of $\Zpnx$, $r\ge 0$, and
$f_1,\ldots,f_r\in \Zpnx $ monic polynomials such that
the residue class of $f_i$ in $\Zpx$ is a power of a 
monic irreducible polynomial $g_i\in \Zpx$,
and $g_1,\ldots,g_r$ are distinct. 

Moreover,
$k\in\natn_0$ is unique, $u\in \Zpnx$ is unique modulo $p^{n-k}\Zpnx$,
and also the $f_i$ are unique (up to ordering) modulo $p^{n-k}\Zpnx$.
\endprofess

\proof
Follows from \cancelpNum, \factoringregularfNum\ and \uniqueprimaryNum.
\endproof

For the the monoid of non-zerodivisors of $\Zpnx$ this means:

\profess{\regularmonoid}
Let $M$ be the multiplicative monoid of non-zerodivisors of
$\Zpnx$, and $U$ its group of units. Let $M'$ be the submonoid
of $M$ consisting of all monic polynomials. Then
$$M\isom U \oplus M'.$$
For every monic irreducible 
polynomial $f\in\Zpx$ let $M_f$ be the submonoid of $M'$
consisting of those monic polynomials $g\in\Zpnx$ whose image
under projection to $\Zpx$ is a power of $f$. 
Then
$$M'\isom \sum_f M_f$$
where $f$ ranges through all monic irreducible polynomials of $\Zpx$.
\endprofess

\proof
Follows by specializing \uniquefactorizationNum\ to non-zerodivisors,
which by \reviewZpnxNum\ are exactly the elements of $\Zpnx$ not
divisible by $p$.
\endproof

The example in the introduction shows that the monoid $M_x$
(corresponding to the irreducible polynomial $x\in\Zpx$) has
infinite elasticity (defined in \ElasticityDefNum\ below).
We will show this of each of the direct summands $M_f$ in
\citeregularmonoid.

\subhead{4. Non-uniqueness of factorization and infinite elasticity}

Irreducible elements are defined in an arbitrary monoid in the same
way as we defined them for the multiplicative monoid of a ring
in \IrredPrimeDefNum.
Note that the monoid of non-zerodivisors of $\Zpnx$, and its
submonoid the monoid of all monic polynomials in $\Zpnx$ are cancellative.

\rmprofess{\ElasticityDef}
Let $(M,\cdot)$ be a cancellative monoid.
{\itemenv
\item{\rm (1)}
For $k\ge 2$, let $\rho_k(M)$ be the supremum of all those $m\in\natn$
for which there exists a product of $k$ irreducibles that can also
be expressed as a product of $m$ irreducibles.
\item{\rm (2)}
The {\em elasticity} of $M$ is $\sup_{k\ge 2}(\rho_k(M)/k)$; in
other words, the elasticity is the supremum of the values $m/k$
such that there exists an element of $M$ that
can be expressed both as a product of $k$ irreducibles and as a
product of $m$ irreducibles.
\par
}
\endrmprofess

%

\profess{\irredexample}
Let $f$ be a monic polynomial in $\intz[x]$ which represents an
irreducible polynomial in $\Zpx$. Let $d=\deg(f)$.
Let $n,k\in\natn$ with $0<k<n$, $m\in\natn$ with $\gcd(m,kd)=1$,
and $c\in\intz$ with $p\dividesnot c$.
Then 
$$f(x)^m + cp^k$$ 
represents an irreducible polynomial in $\Zpnx$.
\endprofess

\proof
Suppose otherwise. Then there exist $g,h,r\in \intz[x]$,
with $g,h$ monic, $g$ irreducible in $\Zpnx$, such that
$$f(x)^m + cp^k = g(x)h(x) + p^nr(x)$$ and $0<\deg g<dm$. 
In $\Zpx$, $g$ is a power of $f$, by unique factorization.
Therefore $\deg g= ds$ with $0<s<m$.

Let $\alpha$ be a zero of $g$.
Let $R$ be the ring of algebraic integers in $\ratq[\alpha]$ and $P$
a prime of $R$ above $p$ with ramification index $e$. Let $v$ be
the normalized valuation on $\ratq[\alpha]$ corresponding to $P$.

Since $f(\alpha)^m = p^n r(\alpha) - c p^k$, we have
$v(f(\alpha)) = ke/m$. As $m$ is relatively prime to $k$, 
$m$ divides $e$. This holds for the ramification index $e_i$
of every prime $P_i$ of $R$ above $p$. Let $f_i$ be the inertia
degree of $P_i$. Then $m$ divides
$\sum_i e_if_i = [\ratq[\alpha]: \ratq] = \deg g = ds$. As $m$ is
relatively prime to $d$, $m$ divides $s$, a contradiction to
$0<s<m$.
\endproof

\profess{\Mfelasticity}
Let $n\ge 2$. Let $f$ be a monic irreducible polynomial in $\Zpx$. 
Let $M_f$ be the submonoid of the multiplicative monoid of $\Zpnx$
defined in \citeregularmonoid.

Then the elasticity of $M_f$ is infinite. 
Moreover, $\rho_2(M_f)=\infty$. 
\endprofess

\proof
In the proof, $f$ will, by abuse of notation, denote a monic 
polynomial in $\intz[x]$ which under canonical projection to $\Zpx$ 
maps to the $f$ in the theorem.

Let $q$ be a prime with $q>\max(n-1,\deg f)$. By \citeirredexample,
$f(x)^q+p^{n-1}$ is irreducible in $\Zpnx$, and 
$$\left(f(x)^q+p^{n-1}\right)2 = f(x)^q (f(x)^q + 2p^{n-1})$$
is an example of factorization of a polynomial in $M_f$ into
either $2$ irreducible factors  - on the left - or more than $q$
irreducible factors - on the right (actually, the number of
irreducible factors on the right is $q+1$ for $p\ne 2$, by
\citeirredexample, and $2q$ for $p=2$.) As $q$ can be made arbitrarily 
large, $\rho_2(M_f)=\infty$, and the elasticity of $M_f$ is infinite.
\endproof

\profess{\elasticity}
Let $M$ be the multiplicative monoid of non-zerodivisors of $\Zpnx$
and $M'$ the submoinoid of monic polynomials. Then the elasticity of
$M'$ is infinite, and $\rho_2(M')=\infty$. Therefore the elasticity
of $M$ is infinite and also $\rho_2(M)=\infty$.
\endprofess

\proof
Indeed, by \citeregularmonoid\ and \citeMfelasticity, $M'$ is an
infinite direct sum of monoids $M_f$ of
infinite elasticity and satisfying $\rho_2(M_f)=\infty$.
\endproof
\bigskip
\goodbreak

\rmprofess{Acknowledgment}
The authors wish to thank F.~Halter-Koch for useful comments.
\endrmprofess